\documentclass[11pt]{amsart}
\usepackage[margin=1in]{geometry}
\usepackage{amsthm}
\usepackage{times,eucal,fullpage,times,amsmath,amsthm,amssymb,mathrsfs,stmaryrd,enumerate,accents}
\usepackage[all]{xy}
\usepackage[usenames,dvipsnames]{color}
\usepackage[colorlinks,citecolor=blue,hyperindex]{hyperref}
\usepackage{epigraph}
\usepackage{graphicx}
\usepackage{fancyhdr} % used by RCS
\usepackage{enumitem}
\newcommand{\be}{\begin{equation}}
\newcommand{\ee}{\end{equation}}
\newcommand{\bes}{\begin{equation*}}
\newcommand{\ees}{\end{equation*}}
\newcommand{\bea}{\begin{eqnarray}}
\newcommand{\eea}{\end{eqnarray}}
\newcommand{\beas}{\begin{eqnarray}}
\newcommand{\eeas}{\end{eqnarray}}
\newcommand{\ben}{\begin{note}}
\newcommand{\een}{\end{note}}
\newcommand{\bexl}{\vskip0.1em\noindent\hrulefill\vskip1em\begin{ExerciseList}}
\newcommand{\eexl}{\end{ExerciseList}\hrulefill}

\newcommand{\bthm}{\begin{theorem}}
\newcommand{\ethm}{\end{theorem}}
\newcommand{\bpro}{\begin{prop}}
\newcommand{\epro}{\end{prop}}
\newcommand{\bcor}{\begin{corollary}}
\newcommand{\ecor}{\end{corollary}}
\newcommand{\bcon}{\begin{conjecture}}
\newcommand{\econ}{\end{conjecture}}
\newcommand{\bp}{\begin{proof}}
\newcommand{\ep}{\end{proof}}
\newcommand{\blem}{\begin{lemma}}
\newcommand{\elem}{\end{lemma}}
\newcommand{\bn}{\begin{note}}
\newcommand{\en}{\end{note}}
\newcommand{\benum}{\begin{enumerate}}
\newcommand{\eenum}{\end{enumerate}}
\newcommand{\bed}{\begin{defn}}
\newcommand{\eed}{\end{defn}}
\newcommand{\brem}{\begin{remark}}
\newcommand{\erem}{\end{remark}}

\newcommand{\btik}{\begin{tikzpicture}\begin{axis}[scale=0.5,axis y line=center, axis x line=middle]}
\newcommand{\etik}{\end{axis}\end{tikzpicture}}

\let\mapsto=\longmapsto

\newcommand{\upperRomannumeral}[1]{\uppercase\expandafter{\romannumeral#1}}

\newtheorem{theorem}[equation]{Theorem}      % (If you want theorem numbered
\newtheorem{lemma}[equation]{Lemma}          %
\newtheorem{corollary}[equation]{Corollary}  %       goes for lemmas, etc.)
\newtheorem{proposition}[equation]{Proposition}

\theoremstyle{definition}

\theoremstyle{definition}
\newtheorem{defn}[equation]{Definition}
\theoremstyle{remark}

\theoremstyle{definition}
\newtheorem{remark}[equation]{Remark}

\numberwithin{equation}{section}
\let\isom=\simeq

\let\tensor=\otimes

\newcommand{\C}{{\mathbb C}}

\newcommand{\F}{{\mathbb F}}

\newcommand{\gal}{{\rm Gal}}

\newcommand{\Q}{{\mathbb Q}}

\newcommand{\tL}{{\tilde{L}}}

\newcommand{\Z}{{\mathbb Z}}

\renewcommand{\int}{\operatorname{int}}
\renewcommand{\O}{{\mathcal O}}

\renewcommand{\wp}{{\mathfrak p}}

\newcommand{\fq}{{\mathbb{F}_q}}

\renewcommand{\bpro}{\begin{proposition}}
\renewcommand{\epro}{\end{proposition}}
\long\def\comment#1\endcomment{}

\begin{document}

\title[]{A method for construction of rational points over elliptic curves}%
\author{Kirti Joshi}%
\address{Math. department, University of Arizona, 617 N Santa Rita, Tucson
85721-0089, USA.} \email{kirti@math.arizona.edu}
%\date{Preliminary Version: \today}

\thanks{}%

%\date{}%
%\dedicatory{}%
%\commby{}%

% ----------------------------------------------------------------
\begin{abstract}
I provide a systematic construction of points (defined over a large number fields) on the Legendre curve over $\Q$: for any odd integer $n\geq 3$ my method constructs $n$ points on the Legendre curve and I show that rank of the subgroup of the Mordell-Weil group they generate is $n$ if $n\geq 7$. I also show that every elliptic curve over any number field admits similar type of points after a finite base extension.
\end{abstract}
\maketitle
%%% EPIGRAPH -----------------------------------------------------
\epigraph{\textit{One pair out of a flock of geese remained sporting in water, and  seven times half the square-root proceeded to the shore, \\tired of the diversion. Tell me, L\={\i}l\={a}vat\={\i}, the number of the flock.}}{from the \textit{L\={\i}l\={a}vat\={\i}} of Bh\={a}skara~\cite{colebrooke}}

%-----------------------------------------------------------------
\tableofcontents

\renewcommand{\wp}{\mathfrak{p}}
\newcommand{\wq}{\mathfrak{q}}
\newcommand{\ok}{\O_K}

% ----------------------------------------------------------------
\section{Introduction}
Recently\footnote{on March 1, 2017} Douglas Ulmer gave a colloquium talk at the University of Arizona on his remarkable construction of explicit rational points on the Legendre elliptic curve over $\fq(t)$. Ulmer's work immediately raises the question of whether such a result exists for elliptic curves over number fields.
This note, and \cite{joshi17-points} in which I provide other methods of constructing elliptic and hyperelliptic curves with interesting rational points, grew out of my attempts to understand Ulmer's construction (see \cite{ulmer14}). I describe here a method which constructs, for any odd integer $n\geq 3$,  $n$ explicit points on the Legendre elliptic
\be 
Y^2=X(X+1)(X+\lambda)
\ee
curve over $\Q$. These points are  defined over   number fields which are also described in my construction. For $n\geq 7$, I show in Theorem~\ref{th:legendre0} that these points span a subgroup of the Mordell-Weil group of rank $n$. In Section~\ref{su:numerical-examples}, I provide numerical examples which illustrate my construction by SAGE computations. My method of construction is  universal (see Theorem~\ref{th:univ}) so it also constructs points on a Legendre elliptic curve over finite fields of characteristics not equal to $2,3$) and quite flexible and many variants are possible (see Remark~\ref{re:variants}). I also show in particular that any elliptic curve over any number field, after a finite extension, carries a similar family of points over similarly defined extensions of the base field (see Theorem~\ref{th:legendre2}) %and describe an even simpler method of constructing similar points on the Legendre curve over $\C(\lambda)$ over a family of finite extensions of $\C(\lambda)$ (see \ref{th:complex}). 
My method of construction uses a remarkable  result of \cite{osada87} whose applicability to construction of points I discovered in the context of \cite{joshi17-points} (where it occurs quite naturally) and after writing the main results of \cite{joshi17-points} I realized its usefulness to the problem of constructing rational points on the Legendre curve (which is the main result of this paper). This result of \cite{osada87} says $X^n-X-1\in\Q[X]$ has Galois group $S_n$ for $n\geq 1$. That this polynomial is irreducible for all $n\geq 1$ was in fact discovered in \cite{selmer56}. I realized during the course of this work that this polynomial also has other remarkable properties: all its roots are global units (see Lemma~\ref{le:units} for a more precise statement).

It is a pleasure to thank Doug Ulmer. This paper was inspired by his stimulating and insightful colloquium. Doug also provided, at very short notice, a number of detailed comments and suggestions which have vastly improved this manuscript (and \cite{joshi17-points}). I also thank Sean Howe for suggesting that a fuller discussion of the action of the symmetric group (in the proof of Theorem~\ref{th:legendre0}) was needed in an earlier version of this manuscript. I thank Klaus Lux for conversations.

\section{Points on the Legendre Elliptic curve over $\Q$}\label{se:legendre}
By a \emph{Legendre elliptic curve over $\Q$}, I mean the curve  $E_\lambda$ defined by the equation:
\be 
E_\lambda:Y^2=X(X+1)(X+\lambda).
\ee
with $\lambda\neq0,1$ and $\lambda\in\Q$. This is a quadratic twist of the curve $Y^2=X(X-1)(X-\lambda)$ which is usually referred to as the Legendre curve in elliptic curve literature.

\bthm\label{th:legendre0}
Consider the Legendre curve over $\Q$ (for $\lambda\neq 0,1$ in $\Q$):
$$E_\lambda:Y^2=X(X+1)(X+\lambda),$$
and let $O$ denote the point at infinity on $E_\lambda$.
Let $L_n/\Q$ be the splitting field of $f(X)=X^n-X-1$ for odd $n\geq 5$ and let  $u\in L$ be any  root of $f(X)=0$ in $L_n$. Let $\tL_n$ be the smallest extension  of $L_n$ which contains $\sqrt{(u+\lambda)}$ for all roots $u$ of $f(X)$.  Then 
\benum[label={\bf(\arabic{*})}]
\item\label{th:points-1} For each root $u$ of $f(X)$, the point $P=(u,u^{(n+1)/2}\sqrt{(u+\lambda)})$ is an $\tL_n$-rational point on $E_\lambda$.
\item\label{th:points-2} For all $n\geq 7$, these $n$ points generate a subgroup of  $E_\lambda(\tL_n)$ of rank equal to $n$.
\eenum
\ethm
\bp 
Since $u^n-u-1=0$ for any root $u$ of $f(X)$  one has $u^n=u+1$. Thus the claim that $P$ is a $\tL_n$-rational point on $E_\lambda$ follows from:
\beas 
u^n&=&u+1\\
(u^{(n+1)/2}\sqrt{u+\lambda})^2&=&u\cdot (u+1) (u+\lambda).
\eeas
This proves \ref{th:points-1}.

Now observe that the Galois group of $f(X)=X^n-X-1\in \Q[X]$ is $S_n$ by \cite{osada87}, further by \cite[Theorem 1]{selmer56} this polynomial is irreducible for all $n$ so the Galois action on the roots is transitive. First of all let me point out that for any root $u$ of $X^n-X-1$ and any $\lambda\in\Q$ the extension $\Q(u)[\sqrt{u+\lambda}]/\Q(u)$ has degree equal to two over $\Q(u)$. Indeed, suppose this is not the case then $\sqrt{u+\lambda}\in\Q(u)$ and hence one has
$$\sqrt{u+\lambda}=a_0+a_1u+\cdots+a_{n-1}u^{n-1},$$
and squaring this and using the fact that $u^n=u+1$ one sees that $u$ satisfies an equation of degree less than $n$ over $\Q$ (for example the highest exponent on the right is $u^{2n-2}=u^n\cdot u^{n-2}=(u+1)\cdot u^{n-2}$). But this is impossible as $X^n-X-1$ is irreducible  (by \cite{selmer56}) and hence this polynomial is the minimal polynomial of $u$ over $\Q$.

Further note that $\tL_n/\Q$ is Galois: the polynomial $\prod_{i=1}^n(X^2-u_i-\lambda)\in L_n[X]$ is  clearly $\gal(L_n/\Q)$-invariant and hence in fact has coefficients in $\Q$ and its splitting field over $\Q$ is $\tL_n$. Note that any permutation of $u_1,\ldots,u_n$ also provides a permutation of $\sqrt{u_1+\lambda},\cdots,\sqrt{u_n+\lambda}$ and in particular  a field automorphism of $\tL_n/\Q$. So $S_n$ is a subgroup of automorphisms of $\tL_n/\Q$ and the subgroup $N$ generated by the involutions $\sqrt{u_i+\lambda}\to -\sqrt{u_i+\lambda}$ is a normal subgroup as its fixed field  $L_n$ and $L_n/\Q$ is Galois. Moreover $N$ does not intersect the subgroup $S_n$. Indeed $S_n\cap N$ is normal in $S_n$ and as $n\geq 7>4$ one sees that $S_n\cap N$ is either $S_n,A_n$ or $\{1\}$.  Suppose the intersection is non-trivial. As $\tL_n$ is obtained by successively 
attaching the square roots  $\sqrt{u_1+\lambda},\ldots,\sqrt{u_n+\lambda}$ to $L_n$ so  $\tL_n/L_n$ is a solvable extension and hence  $N=\gal(\tL_n/L_n)$ is solvable. Thus $N\cap S_n$ is also solvable which is impossible as $n\geq 7>4$ and hence this intersection is trivial. So $\gal(\tL_n/\Q)$ is a semi-direct product of $S_n$ and $N$. In particular one sees that $\tL_n$ is equipped with an action of the quotient group $S_n=\gal(L_n/\Q)$ and this action also provides an action of $S_n$ on the Mordell-Weil group $E(\tL_n)$.

Now let me describe the action of the Galois groups of $\tL_n/L_n$ and $\gal(L_n/\Q)=S_n$  on these points. Let $u_1,\ldots,u_n$ be the $n$ distinct roots of $f(X)=X^n-X-1$ in its splitting field $L_n$ and let $P_i=(u_i,u_i^{(n+1)/2}\sqrt{u_i+\lambda})$ for $1\leq i\leq n$. The automorphisms $\tau_i:\sqrt{u_i+\lambda}\mapsto -\sqrt{u_i+\lambda}$ for $1\leq i\leq n$ of $\tL_n/L_n$ maps the point $P_i=(u_i,u_i^{(n+1)/2}\sqrt{u_i+\lambda})$ to $-P_i\in E_\lambda(\tL_n)$ while and $S_n$ operates transitively by permuting the $P_i$ for $1\leq i\leq n$ (action induced by transitivity of its action on the roots $u_1,\ldots, u_n)$. Hence the action of $S_n$ on the subgroup of $E_\lambda(\tL_n)$ generated by the points $P=(u,u^{(n+1)/2}\sqrt{u+\lambda})$ is transitive.

First suppose that the points $P_1,\ldots,P_n$ are all non-torsion. Consider the vector subspace spanned by $P_1,\ldots,P_n$ in  $E(\tL_n)\tensor_{\Z}\Q$. 
I claim that this subspace of $E(\tL_n)\tensor\Q$ has rank $n$: in fact as $S_n$ acts transitively on $P_1,\ldots,P_n$ the subspace these points generate is the tautological permutation representation of $S_n$ (of dimension $n$). Indeed if these points are non-torsion then the symmetric group $S_n$ action on a non-zero vector space (generated by these points). So one has a representation of the symmetric group on a $\Q$-vector space. This action is faithful. For if not then the action factors through the unique quotient $S_n\to \Z/2$ by the alternating group. But the action of $\Z/2$ on $P_1,\ldots,P_n$ cannot be transitive as $n\geq 7$. So the action is faithful. On the other hand by \cite{dickson08} the dimension of smallest faithful representation of $S_n$ for $n\geq 7$ is $n-1$. So this vector space has rank at least $n-1$. But one can do better: the dimension is $n$.

Indeed, there is a surjection (compatible with action of $S_n$) from the tautological representation of $S_n$ on a $\Q$-vector space of dimension $n$ on which $S_n$ acts by permuting the standard basis, to the subspace of $E(\tL_n)$ generated by these points (the surjection maps the $n$ standard basis vectors to these $n$ points). This tautological permutation representation of $S_n$ is not irreducible: it is the direct sum of a one dimensional trivial representation (generated by the sum of the standard basis vectors) and an irreducible representation of dimension $n-1$ (the standard representation of $S_n$). So the kernel of this surjection is either trivial (in which case the rank is $n$) or of dimension one or dimension $n-1$ and I have already shown that the image has dimension at least $n-1$. So suppose the kernel has dimension one. Then this means that the invariant point $Q=P_1+\cdots+P_n\in E(\tL_n)$ is torsion. I claim that this is not the case. If $Q$ is torsion then by applying the automorphism $\sqrt{u_2+\lambda}\mapsto -\sqrt{u_2+\lambda},\cdots,\sqrt{u_n+\lambda}\mapsto\sqrt{u_n+\lambda}$  of $\tL_n$ which replaces $P_2,\ldots,P_n$ by $-P_2,\ldots,-P_n$ one sees that $Q'=P_1-P_2-\cdots-P_n$ is also torsion and hence $Q+Q'=2P_1$ is torsion; and so $P_1$ is torsion. So using transitivity of the $S_n$ action one sees that all $P_1,\ldots, P_n$ are torsion. This contradicts our assumption that these points are not torsion.

So the mapping is an isomorphism. Thus the points $P_1,\ldots,P_n$ generate the tautological permutation representation of $S_n$ in $E(\tL_n)$ and hence the rank of the subgroup of $E(\tL_n)$ generated by all these points is $n$.

Now let me prove that the points $P_1,\ldots,P_n$ are non-torsion. Suppose one of these points is of finite order. Then by Galois action one sees that all of them must be torsion. Now suppose $NP=O$. Then one can assume $N>2$ as the points are not two-torsion and one can assume that $N$ is minimal such that $NP=O$. Then by transitivity of the action of $S_n$, for all $1\leq i\leq n$, one has $NP_i=O$ and $N$ is minimal with this property. I first claim that $N$ is a power of $2$. Suppose this is not the case. Then there is some odd prime $p$ dividing $N$. Write $N=pM$ with $M\in\Z$ so that $M$ is divisible (if at all) by a smaller power of $p$. Now the point $MP$ is of order exactly $p$, in particular $MP_1,\ldots,MP_n$ are not equal to $O$ and each is of order $p$. I claim that these $n$ points are all distinct. For if, say, $MP_1=MP_2$ then by applying an automorphism of $\tL_n$ which maps $P_2$ to $-P_2$ while keeping all the other points fixed, one gets
$MP_1=-MP_2=MP_2$ so $MP_2=-MP_2$ and hence $2MP_2=O$ and $2M<N$ which contradicts the minimality of $N$. Thus all the points $MP_1,\ldots,MP_n$ are pairwise  distinct. Since these points are of order $p$, they span a subspace of the two dimensional vector space $E[p]$ of the $p$-torsion points of $E$. If this subspace is one dimensional then one gets a homomorphism $S_n\to \F_p^*$ and such a homomorphism factors through the quotient $S_n/A_n\isom\Z/2\to\F_p^*$, or this subspace is two dimensional and one gets a homomorphism from $S_n$ to ${\rm GL}_2(\F_p)$. But by a well-known result of \cite{dickson08}, if $n\geq 7$  there are no faithful homomorphisms from $S_n\to {\rm GL}_2(\Z/p)$ so any such homomorphism is either trivial or factors through the quotient $S_n\to S_n/A_n=\Z/2$ of the symmetric group by the alternating group. Thus in either case action of $S_n$ on these points cannot be transitive.

In either case one has arrived at a contradiction. Thus the only possibility is that this subspace of $E[p]$ is the zero subspace, so the points are of order $M<N$ which again contradicts the minimality of $N$. So $N$ is not divisible by any odd prime. So $N$ is a power of two as claimed.

Now suppose $N=2^m$. Then one has a homomorphism from $S_n$ to $GL_2(\Z/2^m\Z)$. But as the latter  group is solvable while the former is not, so  there are no injective homomorphisms from $S_n$ to this group and hence any such homomorphism factors through the canonical quotient $S_n\to \Z/2\to{\rm GL_2}(\Z/2^m)$.
So again $S_n$ cannot act transitively on these points. So one sees that none of these points $P_1,\ldots,P_n$ are of finite order. 
This completes the proof of \ref{th:points-2}.
\ep

\section{Numerical Examples}\label{su:numerical-examples}
Let me give one or two numerical examples, computed using \cite{sage}. 
\subsection{$\bf{\lambda=5,n=3}$}

Consider 
$$y^2=x(x+1)(x+5)=x^3+6x^2+5x$$ so $\lambda=5$. Let $u$ be a root of $x^3-x-1=0$, let $K=\Q(u)$ and $L=K(t)$ where $t=\sqrt{u+5})$. Then $P=(u,u^{(3+1)/2}\sqrt{u+5})=(u,u^2t)$.
Now $$2P=(141/484*u^2 + 505/484*u - 851/484, (1841/117128*u^2 + 7758/14641*u - 126789/117128)*t )$$
and
\begin{multline}
3P=(6447133488/817674025*u^2 + 395644561/817674025*u - 1700009296/817674025, \\  (-244248318190031/23381388744875*u^2 - 106010387784432/23381388744875*u + \\ 69227442607152/23381388744875)*t).
\end{multline}
Moreover one can compute the torsion subgroup:
$${\rm Torsion}(E(L))=\Z/2\times\Z/2,$$
and as the curve already has all its torsion defined over $\Q$,  
so the point $P$ is non-torsion in $E(L)$. Thus $P$ is also of infinite order  in $L\subset \tL_3$. 

%Moreover using SAGE one also finds that 
%$${\rm Torsion\,}E_\lambda(\tL_3)=\Z/2\times\Z/2$$
%and one can also compute the height pairing matrix of the three constructed points and it turns out to be
%$$
%\begin{pmatrix}
%	     1.00342559193039 &    0.000000000000000 & 1.11022302462516\times 10^{-16} \\
%	   0.000000000000000   &   1.00342559193039 & -1.11022302462516\times 10^{-16} \\
%	 1.11022302462516\times 10^{-16} & -1.11022302462516\times 10^{-16} &     1.00342559193039
%\end{pmatrix}
%$$
%and its determinant is $1.01031202002960$. Thus one also finds computationally that $P_1,P_2,P_3$ are linearly independent in the Mordell-Weil group.
\subsection{$\bf{\lambda=7,n=5}$}

Now consider $$y^2=x(x+1)(x+7)=x^3+8x^2+7x,$$ and let $u$ be a root of $x^5-x-1=0$ and $K=\Q(u)$, $L=K(\sqrt{u+7})=K(t)$. Then the point $P=(u,u^3*t)$ and 
$${\rm Torsion}(E(L))=\Z/2\times\Z/2$$
So again $P$ is of infinite order in $E(L)$ and hence in $E(\tL_5)$.
One can compute
\begin{multline}
2P = (16843/67204*u^4 + 25128/16801*u^3 - 24687/16801*u^2 + 103201/67204*u - 151215/67204, \\ 
(-406110339/564547202*u^4 + 39907057/2258188808*u^3 + \\
 423990243/564547202*u^2 - 3681558523/2258188808*u + 421557246/282273601)*t)
\end{multline}

\section{Construction of points is universal}
My construction is in fact universal:
\bthm\label{th:univ}
For every odd integer $n\geq 3$ there exists  a surjective homomorphism of $\Z$-algebras:
$$ 
\frac{\Z[\lambda][X,Y]}{(Y^2-X(X+1)(X+\lambda))}\to \frac{\Z[\lambda,U,T]}{(U^n-U-1,T^2-(U+\lambda))}
$$
given by 
\beas 
\lambda &\mapsto & \lambda \\
X & \mapsto & U\\
Y & \mapsto & U^{(n+1)/2}\cdot T=U^{(n+1)/2}\cdot \sqrt{U+\lambda}.
\eeas
\ethm
\bp 
Consider the homomorphism from $\Z[\lambda,X,Y]\to \frac{\Z[\lambda,U,T]}{(U^n-U-1,T^2-(U+\lambda))}$ given by the above formulae. Then this factors through the principal ideal $(Y^2-X(X+1)(X+\lambda))$ because of the identity:
$$(U^{(n+1)/2}T)^2=(U^{(n+1)/2}\sqrt{U+\lambda})^2=U(U+1)(U+\lambda).$$
The only point which remains to be checked is the surjectivity for which it is enough to prove that the image of $U$ in $\frac{\Z[\lambda,U,T]}{(U^n-U-1,T^2-(U+\lambda))}$ is a unit. Since this fact will also be of use later, I record its proof in the following Lemma.
\ep 

\blem\label{le:units}
For all integers $n\geq2$, for any root $u$ of the polynomial $X^n-X-1\in\Q[X]$, $u$ and $u+1,u-1$ are  units in the ring of integers $\O_{L_n}$ of its splitting field $L_n/\Q$. More precisely one has the following:
\benum[label={{\bf(\arabic{*})}}]
\item $u(u^{n-1}-1)=1$,
\item $(u+1)(u^{n-1}-u^{n-2}+\cdots+u^2-u)=1$, and
\item $(u-1)(u^{n-1}+\cdots+u)=1$.
\eenum
\elem
\bp
In fact $u$ is already a unit in $\Z[X]/(X^n-X-1)\subset\O_{L_n}$: if $\wp$ is a non-zero prime ideal which contains $u$ then $\wp$ contains $u^n-u=1$, so $\wp$ contains $1$. Thus $u$ does not be long to any non-zero prime ideal and clearly $u\neq 0$ so $u$ does not belong to any prime ideal of $\Z[X]/(X^n-X-1)$. So $u$ is a unit. So $u$ is also a unit in $\O_{L_n}$. Since $u^n=u+1$, it follows that $u+1$ is a unit. Similarly if $u-1\in\wp$ for some non-zero prime ideal $\wp$ then so is $(u-1)(u^{n-1}+\cdots+1)=u^n-1=u\in\wp$. So again a contradiction. Hence $u-1$ is also a unit.

Now let me prove the more precise version by finding the inverses of $u,u+1,u-1$. Note that $u^n-u-1=0$ says $u^n-u=1$ or $u(u^{n-1}-1)=1$ which is the first equation. Now  the equation for $u-1$ is  verified by using $u^n=u+1$ and:
$$1=(u-1)(u^{n-1}+\cdots+u).$$ Similarly one verifies that
$$(u+1)(u^{n-1}-u^{n-1}+u^{n-2}+\cdots -u)=u^n-u=1.$$
\ep

\begin{remark}
Since it is possible that the results of this paper may also be of interest to readers who may not be familiar with algebraic geometry, let me explain ``universality'' of the construction embodied in Theorem~\ref{th:univ}. Let $R$ be any ring commutative ring with unity. Then the theorem implies in particular that there is a surjection of $R$-algebras
$$ 
\frac{R[\lambda][X,Y]}{(Y^2-X(X+1)(X+\lambda))}\to \frac{R[\lambda,U,T]}{(U^n-U-1,T^2-(U+\lambda))}.
$$
So every elliptic curve $Y^2=X(X+1)(X+\lambda)$,  defined over any commutative ring $R$ with unity is equipped with the point $(u,u^{(n+1)/2}\sqrt{u+\lambda})$ where $u$ is the image of $U$ in $R_n=R[U]/(U^n-U-1)$. This point has coordinates in the ring  $R_n[\sqrt{u+\lambda}]$ (i.e. the smallest ring in which a root $u$ of $U^n-U-1=0$ and the roots of $T^2-(u+\lambda)=0$ exist).
\end{remark}

\begin{remark}
In particular suppose $\fq$ is a finite field with $q$ elements and characteristic $p\neq 2,3$. Then the above construction also produces $n$ points on any Legendre elliptic curve  over $\fq$: $y^2=x(x+1)(x+\lambda)$ with $\lambda\neq 0,1$ in $\fq$. These points live over   specific finite extensions of $\fq$.
For example $x^{15}-x-1$ is irreducible in $\F_{173}[x]$ and hence the elliptic curve $y^2=x(x+1)(x+\lambda)$ (here $\lambda\in\F_{173}$ chosen so that this is an elliptic curve) acquires  the point $(u,u^{8}\sqrt{u+\lambda})$ over the extension $\F_{173}(u,\sqrt{u+\lambda})$ where $u$ is a root of $x^{15}-x-1$ in $\bar{\F}_{173}$.
\end{remark}

\brem\label{re:variants}
There are many variants of my construction. One may use trinomials $X^n-aX-b$ to find rational points on elliptic (and even hyperelliptic curves--hopefully this will be treated in a forthcoming version of this paper). For example \cite{selmer56} asserts that by a result of Serret and Ore, for every prime $p$ and any integer $a$ with $p\nmid a$ the polynomial $x^p-x+a$ is irreducible. So for example taking $p$ to be an odd prime, $a=-\lambda\in\Z$ one gets, for every root $u$ of $f(X)=X^p-X-\lambda$, the points $(u,u^{(p+1)/2}\sqrt{u+1})$ on $y^2=x(x+1)(x+\lambda)$. So many variants of my construction are possible (though surjectivity may hold only over $\Z[1/N]$ for some suitable $N\geq 1$).
\erem
\section{Elliptic Curves over number fields}
The method of construction of points in fact can be applied to any elliptic curve over any number field:
any elliptic curve $E$ over a number field $K$, after passage to the extension $K'=K(E[2])\supseteq K$ of degree at most six, and after applying a suitable automorphism over $K'$, has an equation over $K'$ of the form $Y^2=X(X+A)(X+B)$ for suitable $A,B\in K'$. Hence it suffices to treat elliptic curves with equations of this form. Here is the general statement.
\bthm\label{th:legendre2} 
Let $K$ be a number field and let $E:Y^2=X(X+A)(X+B)$ be any elliptic curve over $K$ with two torsion defined over $K$. Let $f_n(X)=X^n-X-A\in K[X]$ for any odd integer $n\geq 3$ and let $L_n$ be its splitting field. Let $\tL_n$ be the finite extension of $L_n$ generated by $\sqrt{u+B}$ for all the roots $u$ of $f_n(X)$.
Then for every root $u$ of $f_n(X)$, the point $P=(u,\pm u^{(n+1)/2}\sqrt{u+B})$ is an $\tL_n$-rational point on $E$.
\ethm
\bp 
This is clear because:
\beas 
u^n&=&u+A\\
\left(\pm u^{(n+1)/2}\sqrt{u+B}\right)^2&=&u\cdot (u+A) (u+B).
\eeas
%The action of $G_n=\gal(L_n/K)$ on these points is through  its action on the roots of $f_n(X)$ (by construction) and so $G_n$ acts on the $\Q$-subspace generated by these points in the Mordell-Weil group $E(\tL_n)\tensor\Q$ and so the assertion on the rank is immediate. If $f_n(X)$ is irreducible one has a transitive action on the points and so one has the assertion.
\ep

% ----------------------------------------------------------------
\bibliographystyle{plain}
\bibliography{points.bib,../../master/joshi.bib,../../master/master6.bib}
\end{document}